\documentclass[12pt,letterpaper]{amsart}

\usepackage{ae,aecompl}
\usepackage[T1]{fontenc}
\usepackage[latin1]{inputenc}
\usepackage{amsmath}
\usepackage{amsfonts}
\usepackage{amsthm}
\usepackage{mathptmx}
\usepackage{subfigure}
\usepackage{graphicx}
\usepackage{amsrefs}
\usepackage{xspace}

\newcommand{\igstd}[1]{\includegraphics[height=4cm]{#1}}
\newcommand{\igconf}[1]{\includegraphics[height=5cm]{#1}}

\newcommand{\goodgap}{%
  \hspace{\subfigtopskip}%
  \hspace{\subfigbottomskip}}

\newcommand{\Wlog}{w.l.o.g.\xspace}

\makeatletter
\renewcommand{\section}{\@startsection{section}{1}%
  \z@{.7\linespacing\@plus\linespacing}{.5\linespacing}%
  {\normalfont\scshape\bfseries\centering}}
\makeatother

\newtheoremstyle{romthm}%
{}%
{}%
{\itshape}%
{}%
{\scshape}%
{.}%
{ }%
{}%

\newtheoremstyle{romrem}%
{}%
{}%
{}%
{}%
{\itshape}%
{.}%
{ }%
{}%

\theoremstyle{romthm}
\newtheorem{thm}{Theorem}
\newtheorem{lem}{Lemma}

\newtheorem{cor}{Corollary}

\theoremstyle{romrem}
\newtheorem*{rem}{Remark}

\theoremstyle{definition}

\title{There are no realizable $15_4$- and $16_4$-configurations}

\author{J\"urgen Bokowski}
\email{bokowski@mathematik.tu-darmstadt.de}
\author{Lars Schewe}
\email{schewe@mathematik.tu-darmstadt.de}
\address{Technische Universität Darmstadt, Fachbereich Mathematik, Schloßgartenstraße 7, D-64289 Darmstadt}


\begin{document}

\begin{abstract}
  There exist a finite number of natural numbers $n$ for which we do
  not know whether a realizable $n_4$-configuration does exist. We
  settle the two smallest unknown cases $n=15$ and $n=16$. In these
  cases realizable $n_4$-configurations cannot exist even in the more
  general setting of pseudoline-arrangements. The proof in the case
  $n=15$ can be generalized to $n_k$-configurations.  We show that a
  necessary condition for the existence of a realizable
  $n_k$-configuration is that $n>k^2+k-5$ holds.
\end{abstract}
\maketitle

\section{Introduction}
\label{sec:introduction}

Point line configurations have a long history. Levi's
book~\cite{leviCon} about the subject starts with the remark that they
can be considered a starting point for studying combinatorial
geometry. It was also Levi who wrote in 1926 the first known
paper~\cite{leviPG} of pseudoline arrangements, the antecedent of the
general oriented matroid concept. Within this paper we use the latter
concept in the context of configurations. We assume the reader to be
familiar with basic concepts from the theory of oriented matroids in
the rank $3$ case (see for instance~\cite{MR0307027} or
\cite{OM}*{Chapter 6}). For an introduction to the theory of oriented
matroids see also \cite{CUP}.

\begin{figure}[htbp]
  \centering
  \includegraphics[height=5cm]{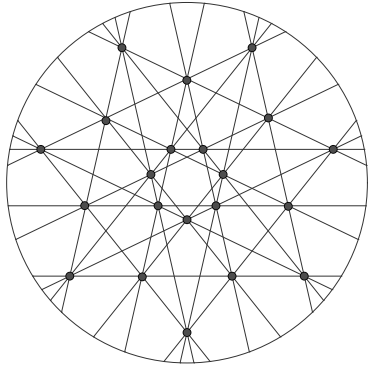}
  \caption{The smallest known realizable $n_4$-configuration}
  \label{fig:214}
\end{figure}

We fix our notation in Section~\ref{sec:proof-theorem-gennk}, however,
an intuitive impression of a realizable $n_4$-configuration can be
obtained from looking at the smallest known example of a realizable
$n_4$-configuration for $n=21$ in Figure~\ref{fig:214}.  A realizable
$n_4$-configuration consists of two $n$ element sets, a set of $n$
points and a set of $n$ lines in the Euclidean plane.  The defining
property of an $n_4$-configuration requires each element of one set to
be incident with precisely 4 elements of the other.

It is known that realizable $n_4$-configurations do not exist for $n
\leq 14$.  For $15 \leq n \leq 20$ and for $n=22, 23, 26, 29, 31, 32,
34, 37, 38, 43$, the existence of realizable $n_4$-configurations is a
long standing problem in this context, whereas for all other n, we do
have realizable configurations (see \cite{MR1913098}). For all even
values that are larger than $21$, realizations with pseudolines are
known. For the values $n=22$ and $n=28$ these can be found see
\cite{MR1067273}. For the other values the realizations were
unpublished; Branko Grünbaum provided us with his drawings (see
Figure~\ref{fig:grun}).

In Section~\ref{sec:proof-theorem-gennk} we prove a general
Theorem~\ref{thm:gennk}, which implies in particular that a realizable
$15_4$-configuration does not exist. The proof for the case $n=16$ in
Section~\ref{sec:case-n=16} is more involved. Our results have been
achieved without using a computer in our final argument, although a
foregoing result of the second author in connection with computations
of Betten and Betten (see \cite{MR1719478}) did use a computer in the
case $n=16$. In Section~\ref{sec:further-results} we give an overview
of this foregoing result.

\section{On general $n_k$-configurations}
\label{sec:proof-theorem-gennk}

An \emph{$n_k$-configuration} (with $k\geq 3$) is a matroid $M$ of
rank 3 on $n$ points such that every line of $M$ has at most $k$
points and each point is contained in exaktly $k$ $k$-point lines. We
say that an $n_k$-configuration is in \emph{general position}, if the
only number of points on a line is either $k$ or $2$.

\begin{thm}\label{thm:gennk}
  A realizable $n_k$-configuration can only exist if $n> k^2+k-5$
  holds. This implies that a realizable $15_4$-configuration does not exist.
\end{thm}

Let $\mathcal{C}$ be a $n_k$-configuration. If an arrangement of
pseudolines has exactly the incidences prescribed by $\mathcal{C}$, we
say that $\mathcal{C}$ is \emph{realizable with pseudo-lines}; or
short \emph{pl-realizable}. By the Folkman-Lawrence representation
theorem (for an elementary proof in the rank-3-case, see
\cite{MR1845486}) this is equivalent to the fact that the matroid $M$
underlying $\mathcal{C}$ is orientable. Our definition of
pl-realizablity implies that every realizable configuration -- that is
one that can be drawn with straight lines in the projective plane --
is also pl-realizable. So pl-realizability is a necessary condition
for realizabilty in the ordinary sense.

We only have to deal with points of $\mathcal{C}$ in general
position:

\begin{rem}
  Let $\mathcal{C}$ be a pl-realizable $n_k$-configuration. Then there
  exists a pl-realizable $n_k$-configuration $\mathcal{C}'$ such that
  $\mathcal{C}'$ is in general position.
\end{rem}

We note that we can always change to polar formulation where we switch
the roles of points and lines.

\begin{proof}[Proof of Theorem~\ref{thm:gennk}]
  The result follows from an application of Euler's formula. We may
  assume that $\mathcal{C}$ is in general position. We assume further
  that we are given a pl-realization of $\mathcal{C}$ on the sphere.
  This induces a graph embedding on the sphere.
  
  We count the number of vertices and edges. The number of vertices is
  given by
  $f_0=2\left(n+\binom{n}{2}-n\binom{k}{2}\right)=n(n-k(k-1)+1)$. The
  number of edges is given by $f_1=2n\left(k + (n-k(k-1)-1)
  \right)=2n(n-k^2+2k-1)$. From Euler's formula we can deduce the
  number of cells: $f_2=f_1-f_0+2$.
  
  A pseudoline-arrangement as a pl-realization of $\mathcal{C}$
  implies that digons are not allowed. By double-counting edge-cell
  incidences we get the following additional inequality:
  $3f_2\leq2f_1$. Plugging in the above expressions for $f_0$, $f_1$,
  and $f_2$ our inequality becomes: \[-{n}^{2}-5\,n+n{k}^{2}+nk+6\leq
  0.\] For fixed $k\geq 3$ and nonnegative $n$ the expression on the
  left-hand-side is monotonically decreasing. For $n=k^2+k-5$ the
  inequality does not hold whereas for $n=k^2+k-4$ the inequality
  holds.
\end{proof}

\begin{rem}
  The proof allows us to replace \emph{realizable} with
  \emph{pl-realizable} in the statement of Theorem~\ref{thm:gennk}.
\end{rem}

\begin{cor}
  Realizable $15_4$-configurations do not exist.
\end{cor}

We are using the fact that two $n_k$-configurations in general
position are $\pi$-equivalent; they have the same Poincaré polynomial
(for a definition, see the book by Orlik and Terao \cite{MR1217488}).

\section{The Case $n=16$}
\label{sec:case-n=16}

The aim of this section is to prove the following theorem:

\begin{thm}\label{thm:164}
  Realizable $16_4$-configurations do not exist.
\end{thm}

We start of with some convenient definitions: We call the intersection
point between pseudolines \emph{crossing}. We call such a point an
\emph{$n$-crossing} if exactly $n$ pseudolines go through that point.
As in the section above we assume that our configuration is in general
position. So, we only have to deal with $2$- and $4$-crossings. We
pick an arbitrary pseudoline as line at infinity in our arrangement
and we denote it with $\infty$.

\begin{figure}[htbp]
  \centering
  \subfigure[]{\igstd{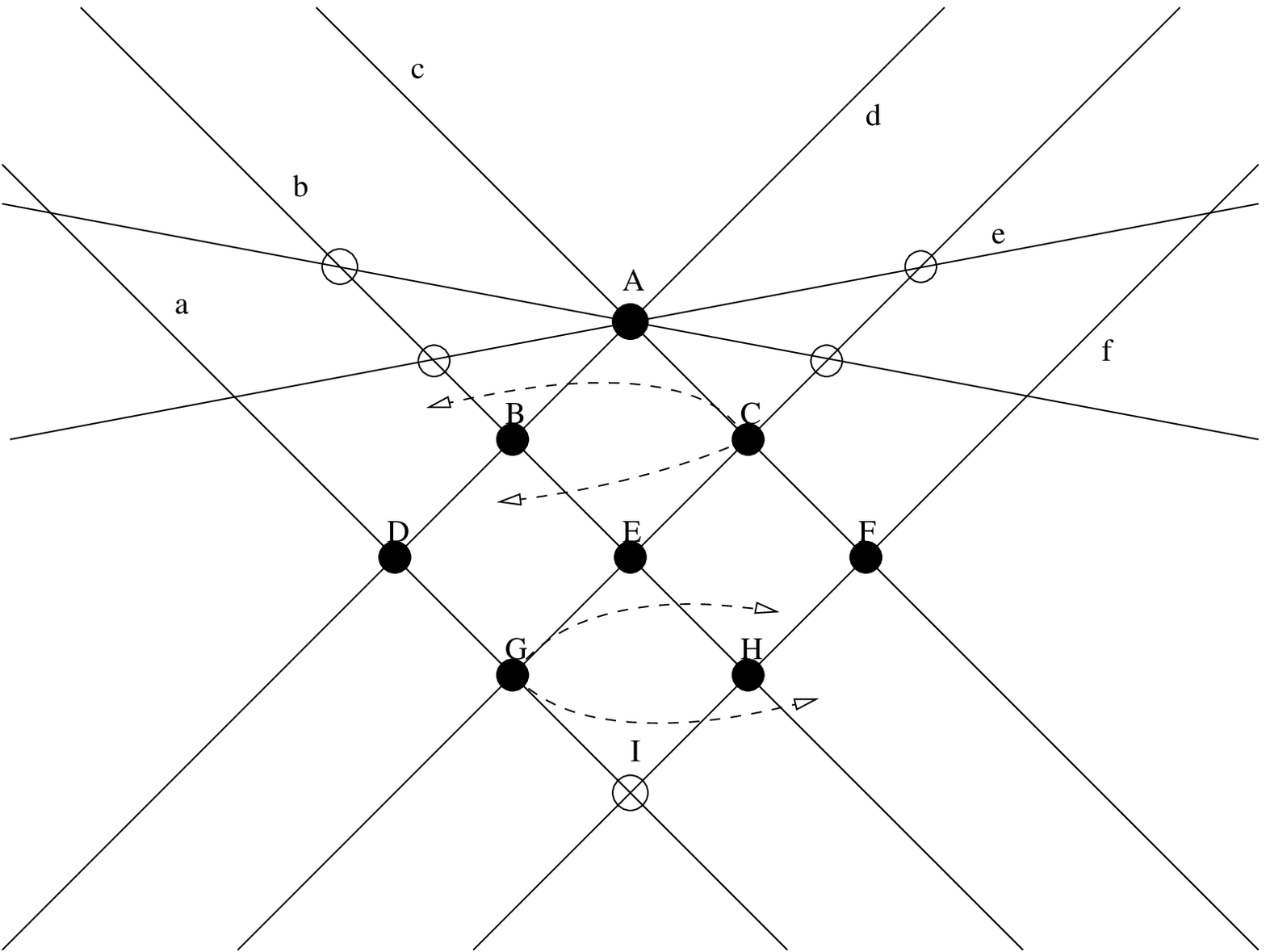}  \label{fig:fall1}}
  \subfigure[]{\igstd{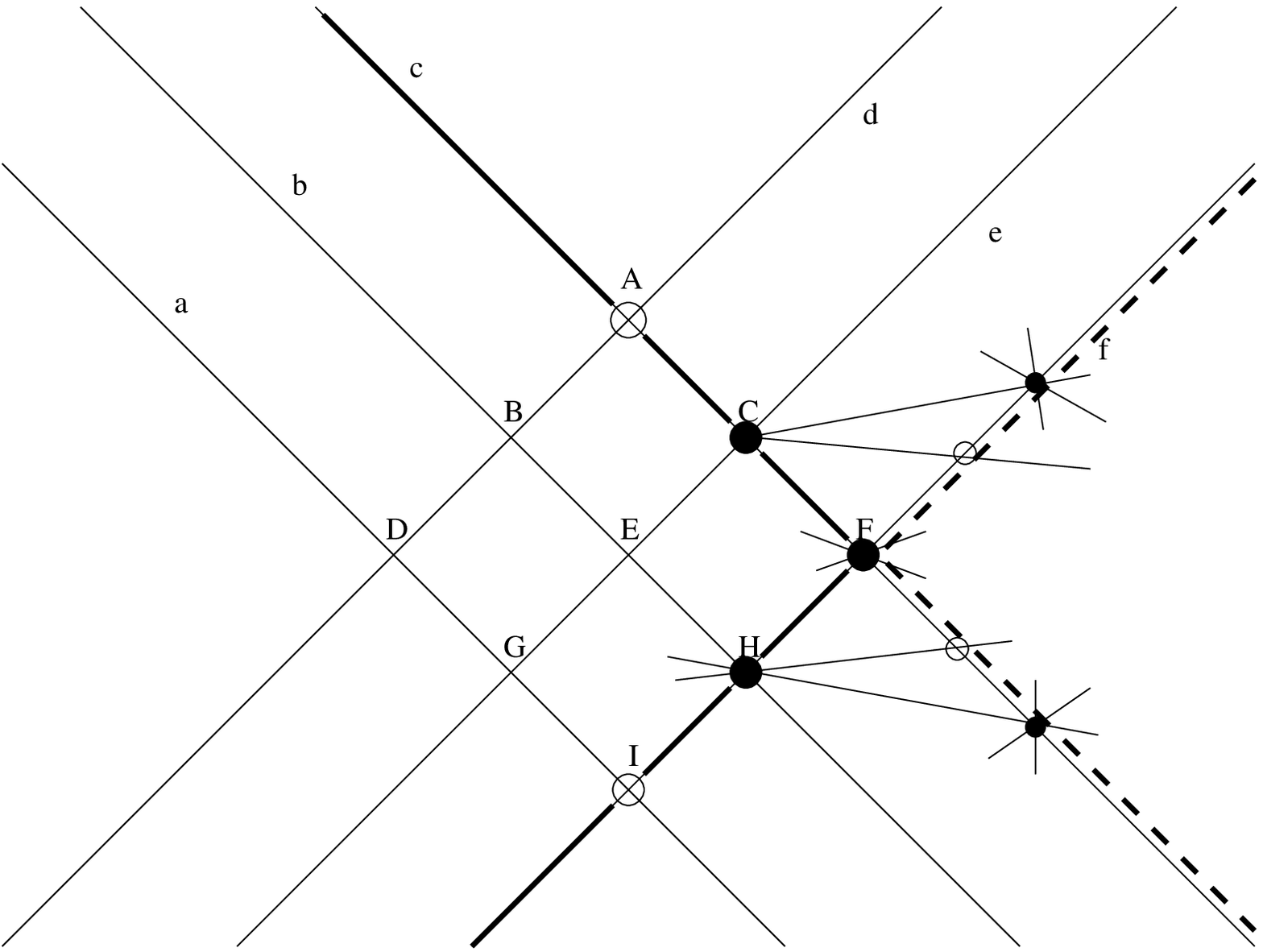}  \label{fig:fall2}}
  \caption{}
\end{figure}

In the case $n=16$ we have to have exactly four $4$-crossings and
three $2$-crossings on each pseudoline. This means that we have two
$4$-crossings that are adjacent on our pseudoline $\infty$.  Now call
the elements that intersect in these two $4$-crossings $a$, $b$, $c$,
and $d$, $e$,$f$, respectively. These six pseudolines have to have
nine additional distinct crossings which we label from $A$ to $I$ as in
Figure~\ref{fig:fall1}. We call these crossings also \emph{grid points}
to distinguish them from other crossings. Note that we have chosen our
starting situation so that no pseudoline in our arrangement can go to
$\infty$ between the pseudolines $c$ and $d$ above $A$ and no
pseudoline can go to $\infty$ between $a$ and $f$ below $I$.

All the pseudolines $a$,\dots,$f$ cross $\infty$ in a 
$4$-crossing. This means that they contain only six further crossings in
total: three $4$-crossings and three $2$-crossings. We further remark
that at most nine and at least six of the grid points have to be
$4$-crossings. 

However, not all nine grid points can be $4$-crossings. This follows
from the following Lemma.

\begin{lem}\label{lem:aisim}
  The grid points $A$ and $I$ cannot be both $4$-crossings. 
\end{lem}
\begin{proof}
  Assume both $A$ and $I$ were $4$-crossings. Then the two non-grid
  pseudolines leaving $A$ have to cross the line $b$ in two distinct
  crossings that both lie above $B$. However, the two
  non-grid pseudolines leaving $I$ have to cross $b$ in two distinct
  crossings as well, but those crossings have to lie below $H$. So we
  get four additional crossings to the crossings $B$, $E$, and $H$ on
  $b$. This is a contradiction.
\end{proof}

If eight grid points are $4$-crossings, we may assume
that one of the points $A$ or $I$ is a $2$-crossing. Now we can show
that the case of eight $4$-crossings in the grid cannot occur.

\begin{lem}\label{lem:amseven}
  At most seven grid points can be $4$-crossings.
\end{lem}
\begin{proof}
  We cannot have nine grid points that are $4$-crossings. This would
  contradict Lemma~\ref{lem:aisim}. So, assume we had eight grid
  points that were $4$-crossings. We may then assume \Wlog that $A$ is
  a $4$-crossing and $I$ is not. So we are in the situation of
  Figure~\ref{fig:fall1}. We can then see that the new pseudolines
  leaving $C$ have to cross $b$ in at least one new crossing. The new
  pseudolines leaving $G$ have to cross $b$ in at least one new
  crossing as well. However, this would give us seven crossings on
  $b$, which is impossible.
\end{proof}

Now we deal with the case that $A$ and $I$ are both
$2$-crossings. 

\begin{lem}\label{lem:aundi}
  The grid points $A$ and $I$ cannot be both $2$-crossings. 
\end{lem}
\begin{proof}
  Assume both were $2$-crossings. Then at most one further grid point
  is a $2$-crossing. So \Wlog we are in the situation of
  Figure~\ref{fig:fall2}. No $4$-crossings can lie in the bold
  $1$-cells of $c$, otherwise we would get too many crossings on
  $f$. By symmetry the same holds for the bold $1$-cells of $g$. So we
  know the $1$-cells in which the further $4$-crossings on $c$
  resp. $g$ lie. We count the number of lines entering the $2$-cell on
  the right which is bordered by the dashed line. This number is ten,
  but only nine lines are leaving this cell to cross $\infty$, which
  is a contradiction.
\end{proof}

From now on we can always assume that $A$ is a $4$-crossing, $I$ is a
$2$-crossing and we have at least one and at most two further grid
points that are $2$-crossings. First we deal with the case that
precisely one further grid point is a $2$-crossing.

\begin{figure}[htbp]
  \centering
  \subfigure[]{\igstd{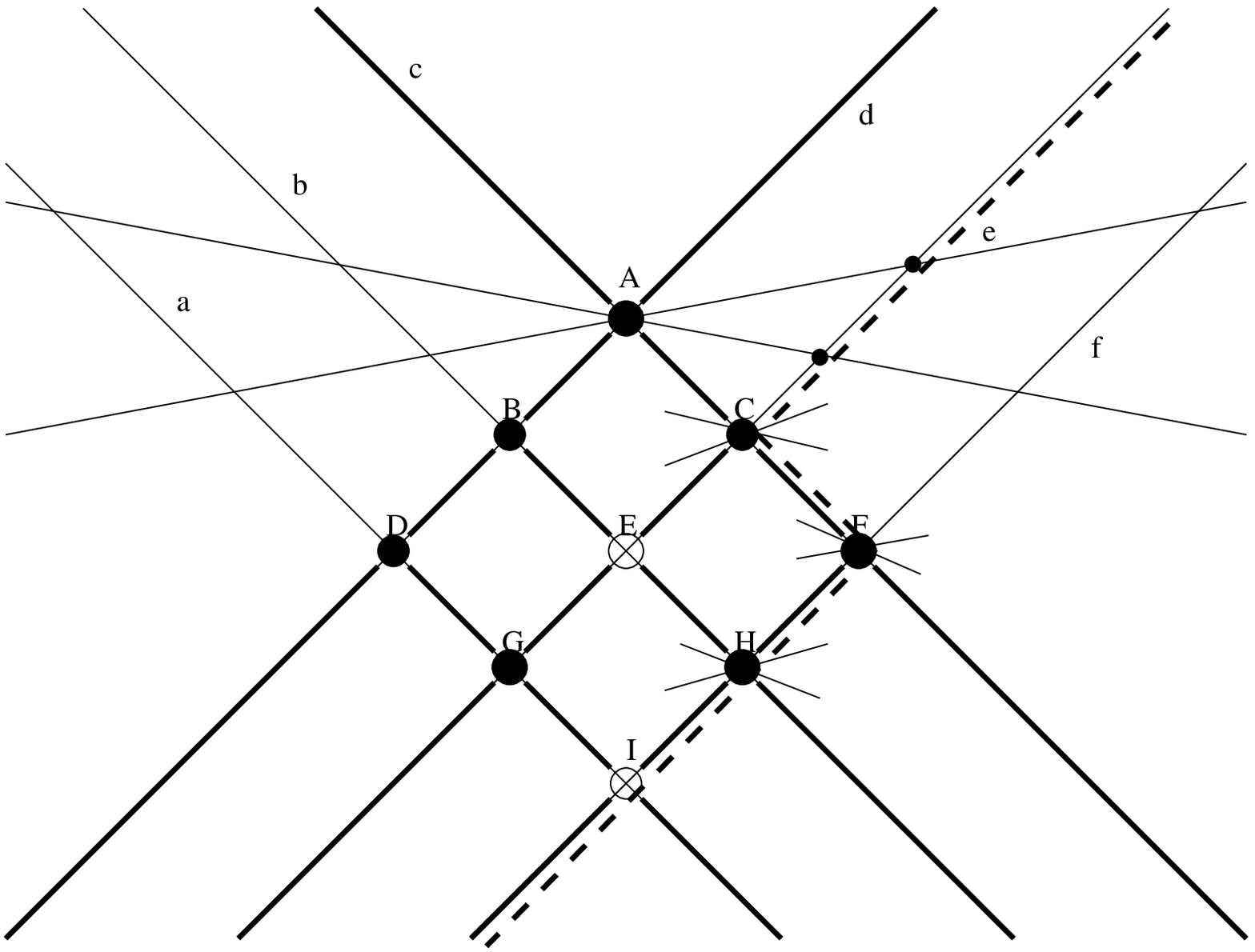}\label{fig:fall3a}}
  \subfigure[]{\igstd{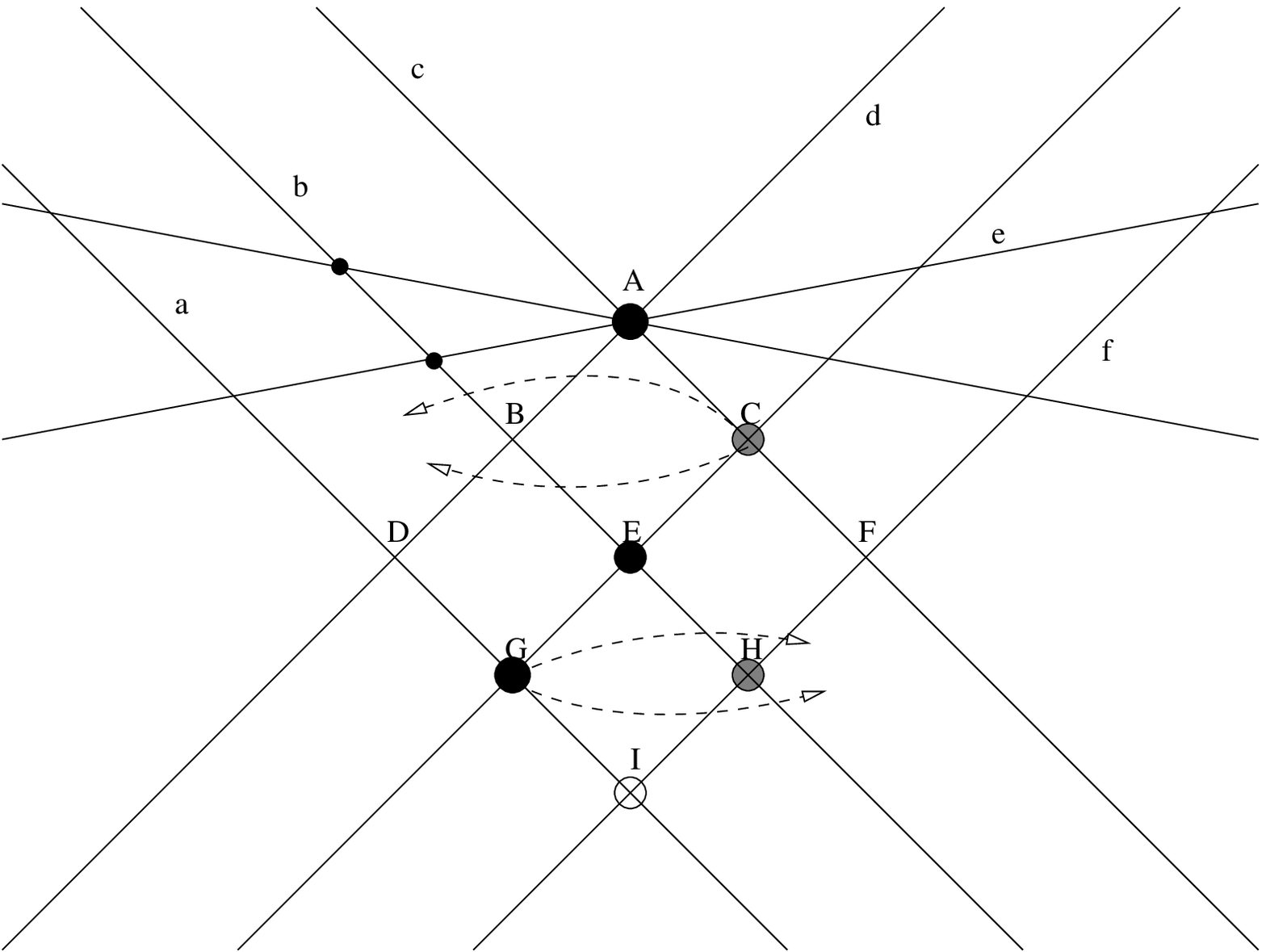}\label{fig:fall3b}}
  \caption{}
\end{figure}

\begin{lem}\label{lem:amsix}
    At most six grid points can be $4$-crossings.
\end{lem}
\begin{proof}
  We deal with two cases seperately: First we assume $E$ is our
  further $2$-crossing. Then we are in the situation of
  Figure~\ref{fig:fall3a}. The bold $1$-cells cannot contain a
  $4$-crossing. This would lead to too many crossings on line $b$, or
  $e$, respectively. So one $4$-crossing has to lie above the point $C$
  on $e$. Now there are again ten pseudolines entering the $2$-cell
  which is bordered by the dashed line, but again only nine
  pseudolines form our arrangement. This is the desired contradiction
  in this case.
  
  So we may assume $E$ is a $4$-crossing. So we are in the situtation
  of Figure~\ref{fig:fall3b}. By symmetry we can assume \Wlog that $G$
  is a $4$-crossing. Now, however, $H$ has to be a $4$-crossing as
  well, otherwise we would get a contradiction. The two lines coming
  from $G$ cross $b$ in two points. If one of these was not $H$, we
  would have eight crossings on $b$. So both $G$ and $H$ have to be
  $4$-crossings.
  
  Now by symmetry we may assume that $C$ is a $4$-crossing. Then the
  two lines coming from $C$ that cross $b$ give at least one new
  crossing on $b$. Together with the crossing that $H$ gave, we have
  eight crossings in total. This is our desired contradiction. So we
  have shown that no seven grid points can be $4$-crossings. Together
  with Lemma~\ref{lem:amseven} we have shown the lemma.
\end{proof}

\begin{figure}[htbp]
  \centering
  \subfigure[]{\igstd{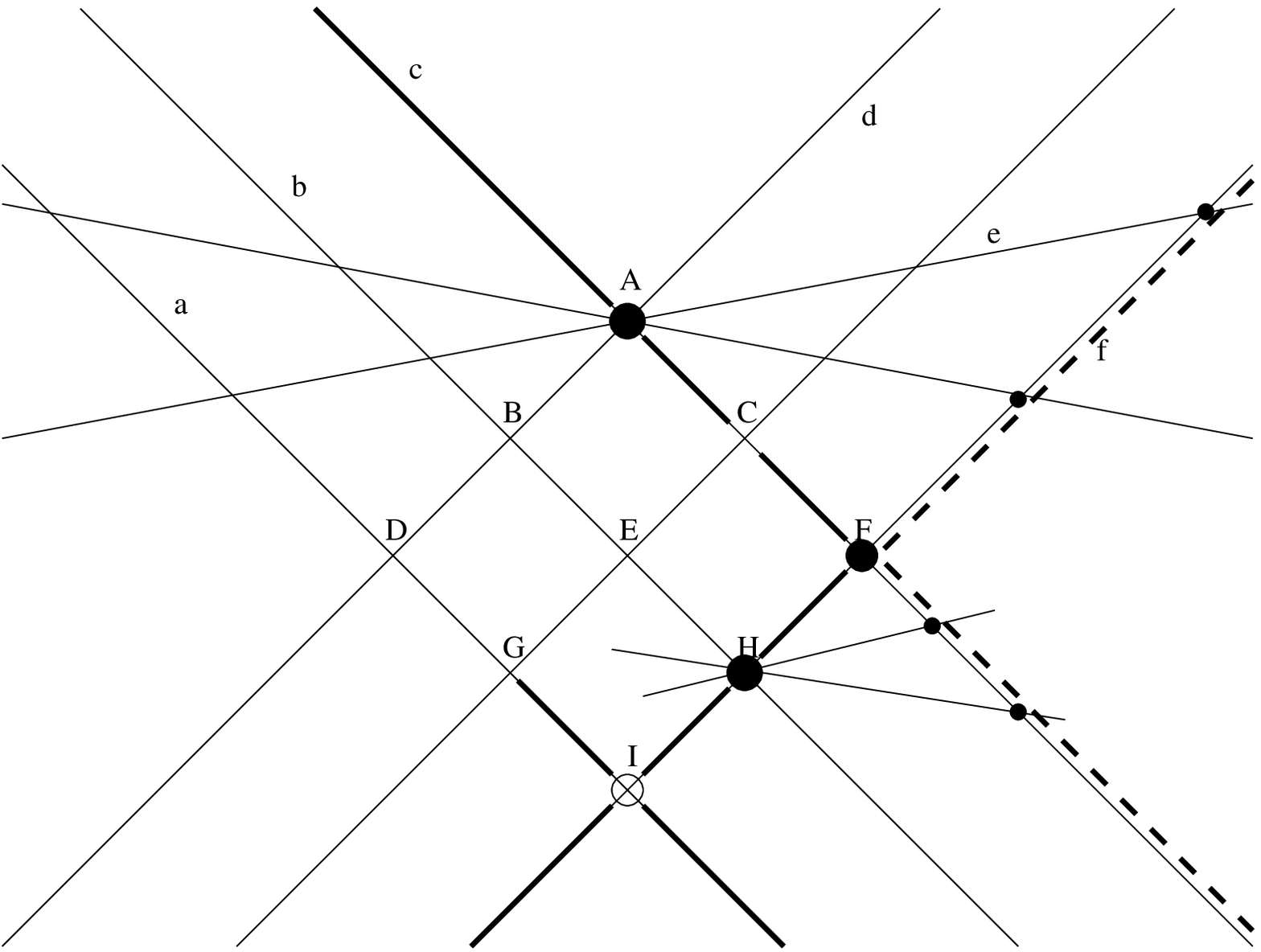}\label{fig:fall4}}
  \subfigure[]{\igstd{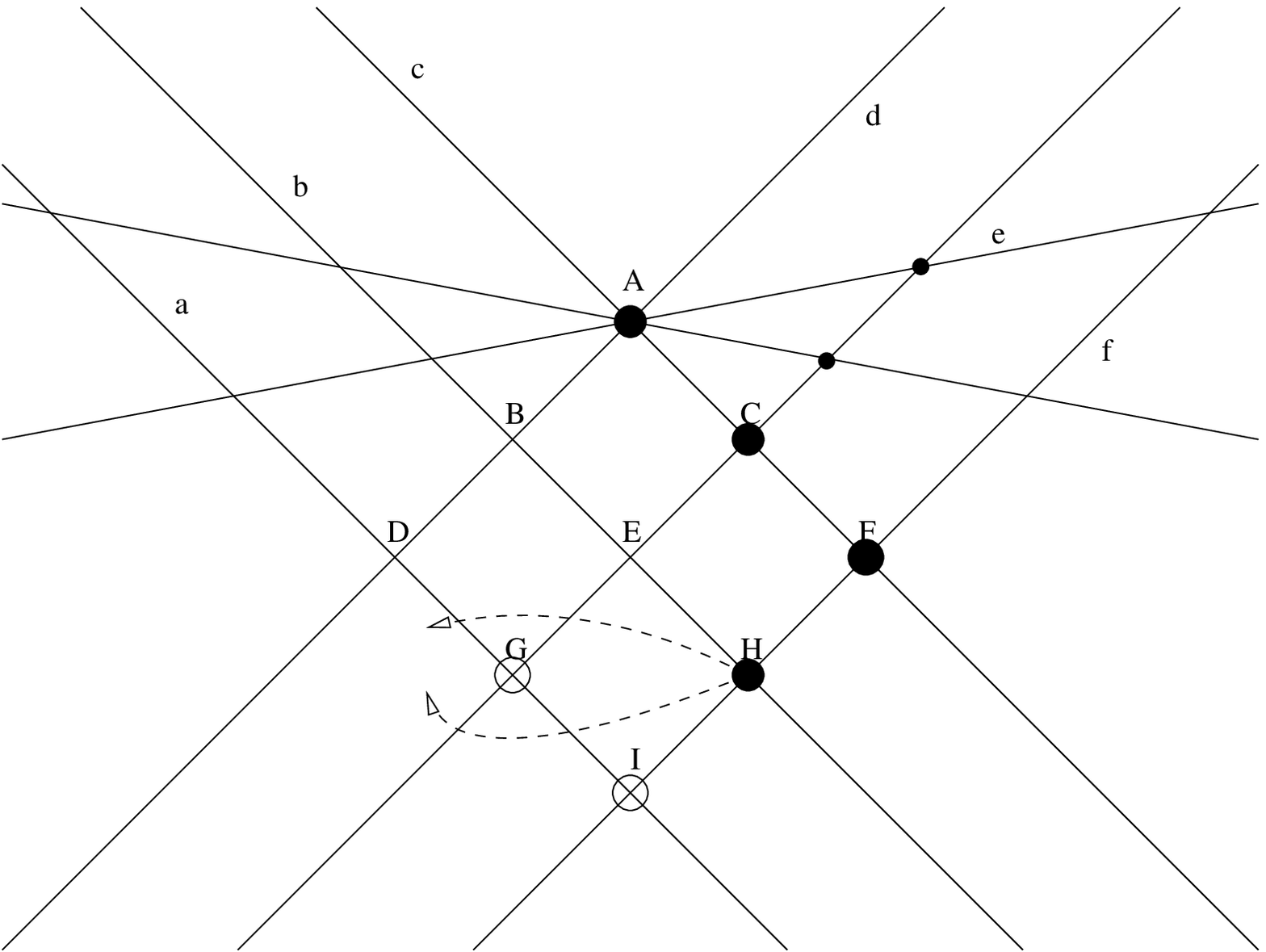}\label{fig:fall5a}}
  \caption{}
\end{figure}

The following lemma reduces the number of possible positions for the
$2$-crossings. 

\begin{lem}\label{lem:afhic}
  Assume $A$, $F$, and $H$ are $4$-crossings and $I$ is a
  $2$-crossing. Then $C$ is a $4$-crossing.
\end{lem}
\begin{proof}
  This is the situation of Figure~\ref{fig:fall4}. Assume $C$ is a
  $2$-crossing. If there was a $4$-crossing on $f$ below $F$, we would
  get too many crossings on $b$. However, no $4$-crossing can lie above
  $F$ on $c$ as well, we would get too many crossings on $f$. So we
  get one $4$-crossing above $F$ on $f$ and one $4$-crossing below $F$
  on $c$. Now the cell bounded by the dashed line is entered by ten
  pseudolines. This is a contradiction. Hence, $C$ has to be a
  $4$-crossing. 
\end{proof}

The next lemma reduces the possibilities further.

\begin{lem}\label{lem:acfhvgi}
  The situation that $A$, $C$, $F$, and $H$ are $4$-crossings and $G$
  and $I$ are $2$-crossings cannot occur.
\end{lem}
\begin{proof}
  This is the situation of Figure~\ref{fig:fall5a}. As can be seen in
  the figure, the lines coming from $H$ that cross $a$ give too many
  crossings on $e$. Hence the situation cannot occur.
\end{proof}

\begin{figure}[htbp]
  \centering
  \subfigure[]{\igstd{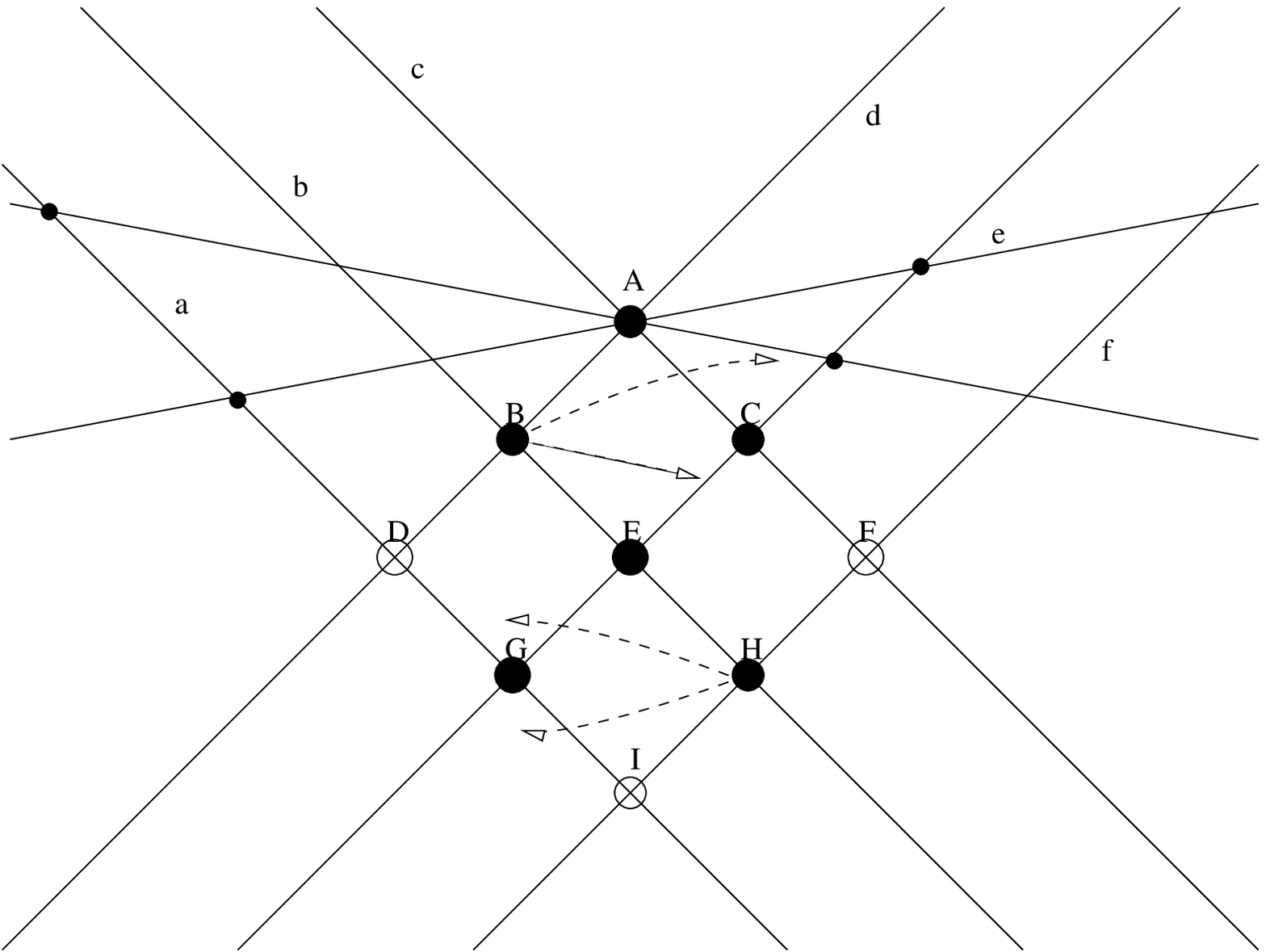}}
  \subfigure[]{\igstd{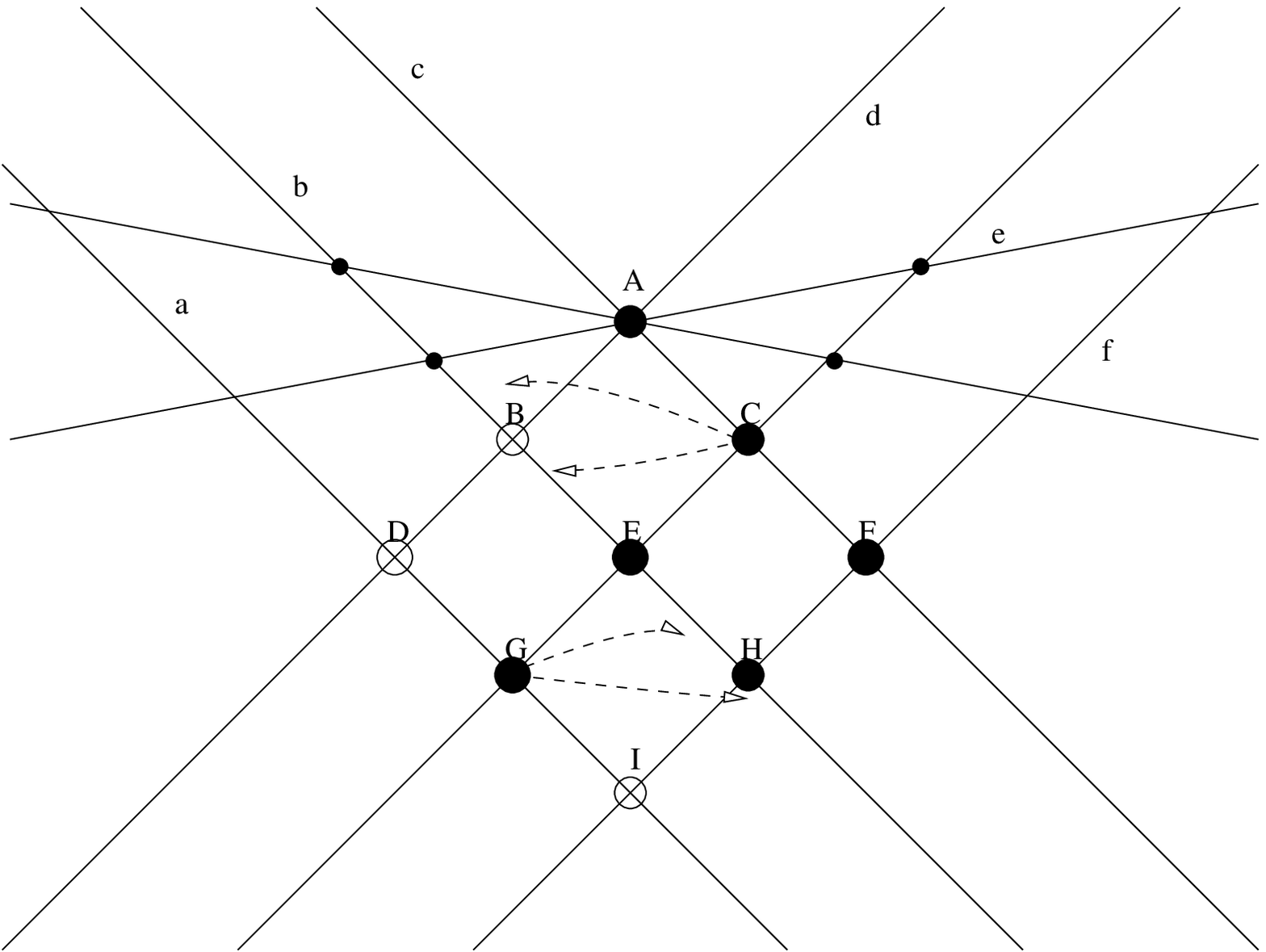}}
  \caption{}
  \label{fig:fall5bc}
\end{figure}

Now only four cases remain, for the first two of them (see
Figure~\ref{fig:fall5bc}) we refer only to the figure. The other two
cases are considerably harder.

\begin{figure}[htbp]
  \centering
  \igstd{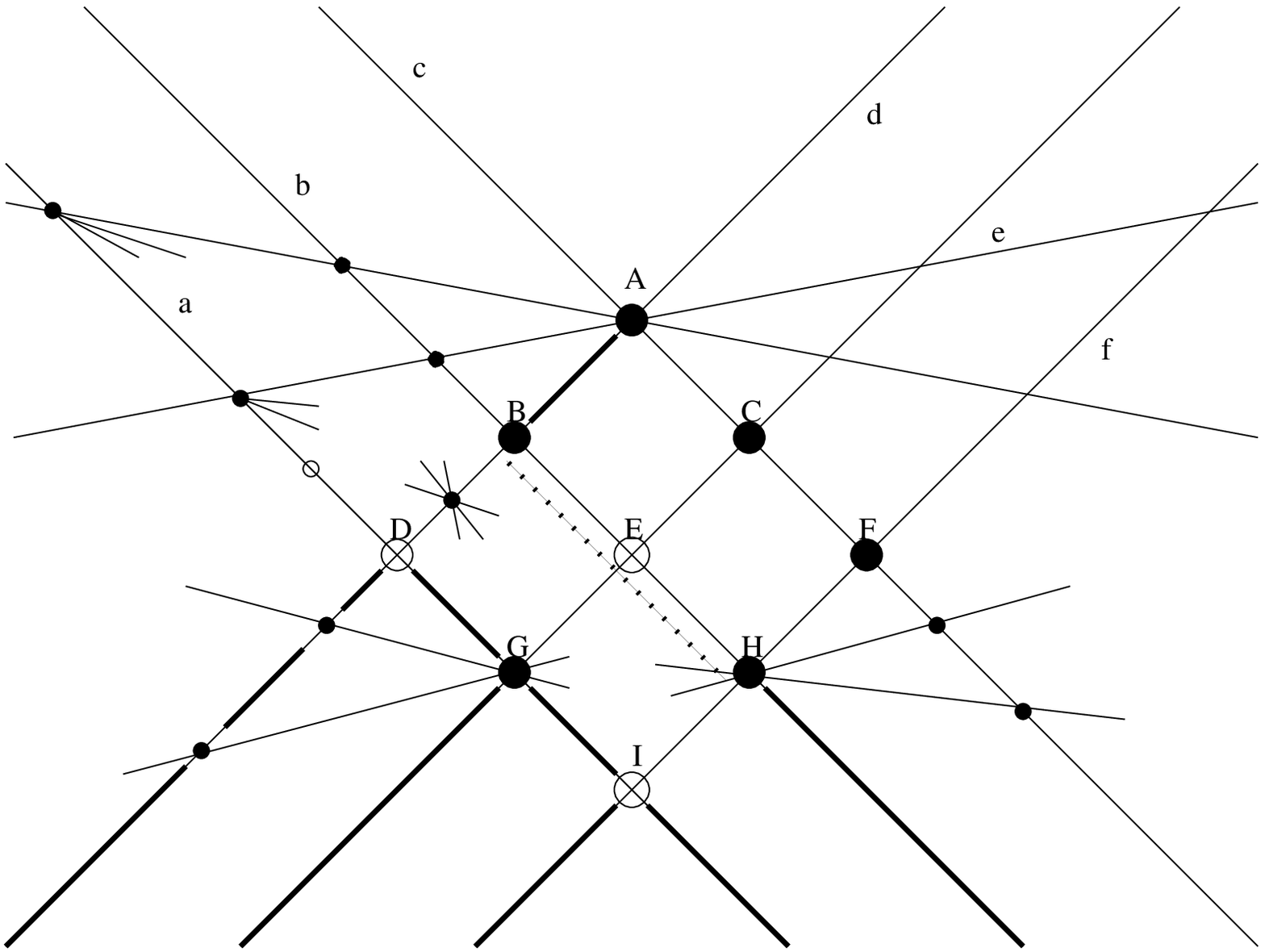}
  \caption{}
  \label{fig:fall6}
\end{figure}

\begin{lem}
  The situation, that the grid points $D$, $E$, and $I$ are
  $2$-crossings cannot occur. 
\end{lem}
\begin{proof}
  This is the situation of Figure~\ref{fig:fall6}. No $4$-crossings
  may lie on $b$ below $H$. Also no $4$-crossings may lie on $d$ above
  $B$. Additionally, no $4$-crossing lies on $a$ below $D$. So the two
  further $4$-crossings on $a$ have to lie above $D$. This means,
  however, that no $4$-crossing on $d$ can lie below $D$. Hence, the
  missing $4$-crossing on $d$ has to lie between $D$ and $B$. Now we
  have five further lines -- coming from the new $4$-crossing and $G$
  -- that need to cross $b$ below $B$. We only have one further
  $4$-crossing that lies on $b$, which can only lie in the segment
  denoted by the dotted line. This means that we only have two
  possible exit points for the above mentioned five lines. This,
  however, leads to a contradiction. The three lines coming from the
  new $4$-crossing on $d$ need to cross $b$ in three pairwise distinct
  crossings.
\end{proof}

Now, we take on the last -- and hardest -- case.

\begin{figure}[htbp]
  \centering
  \subfigure[]{\igstd{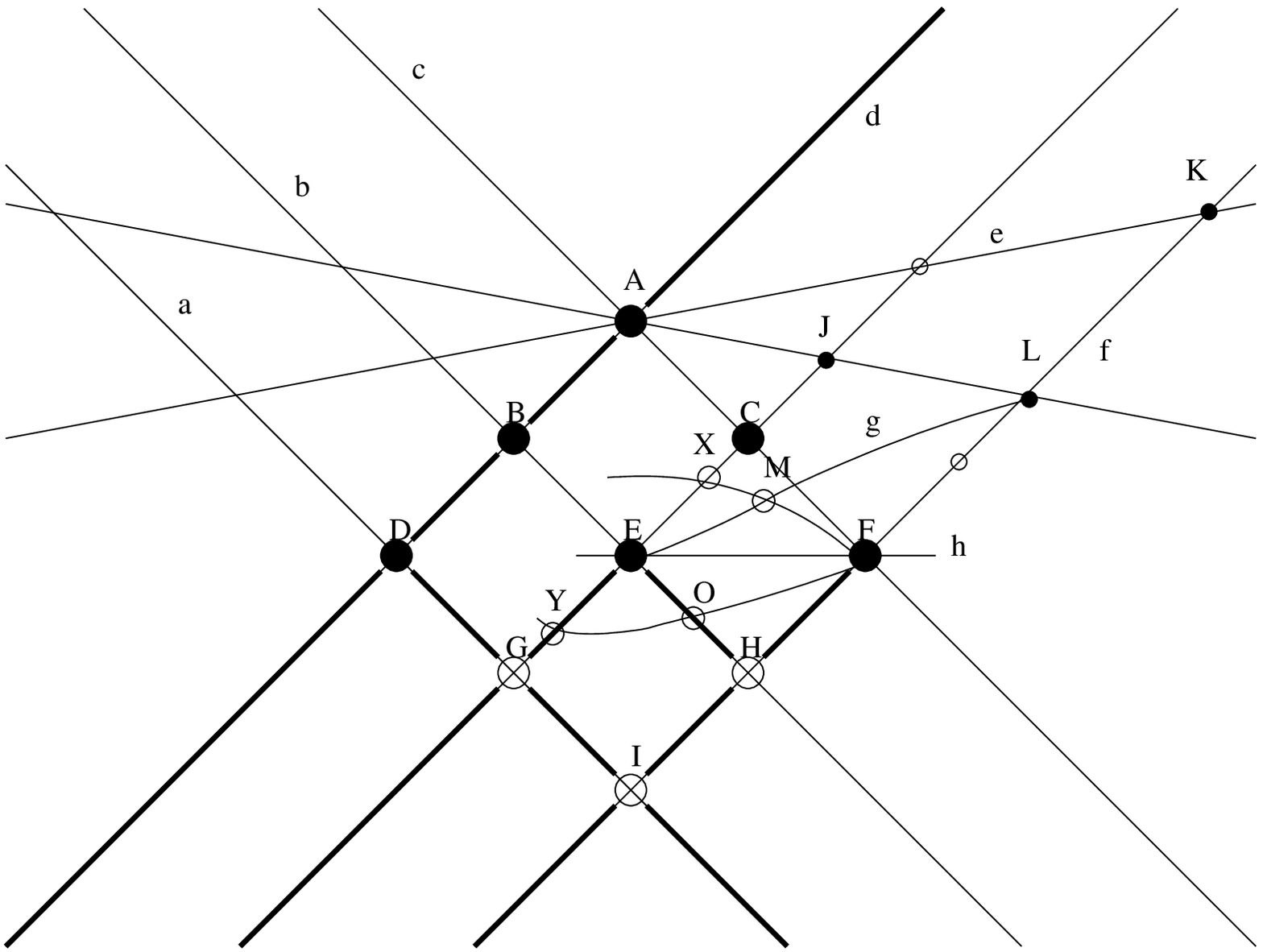}\label{fig:fall7a}}
  \subfigure[]{\igstd{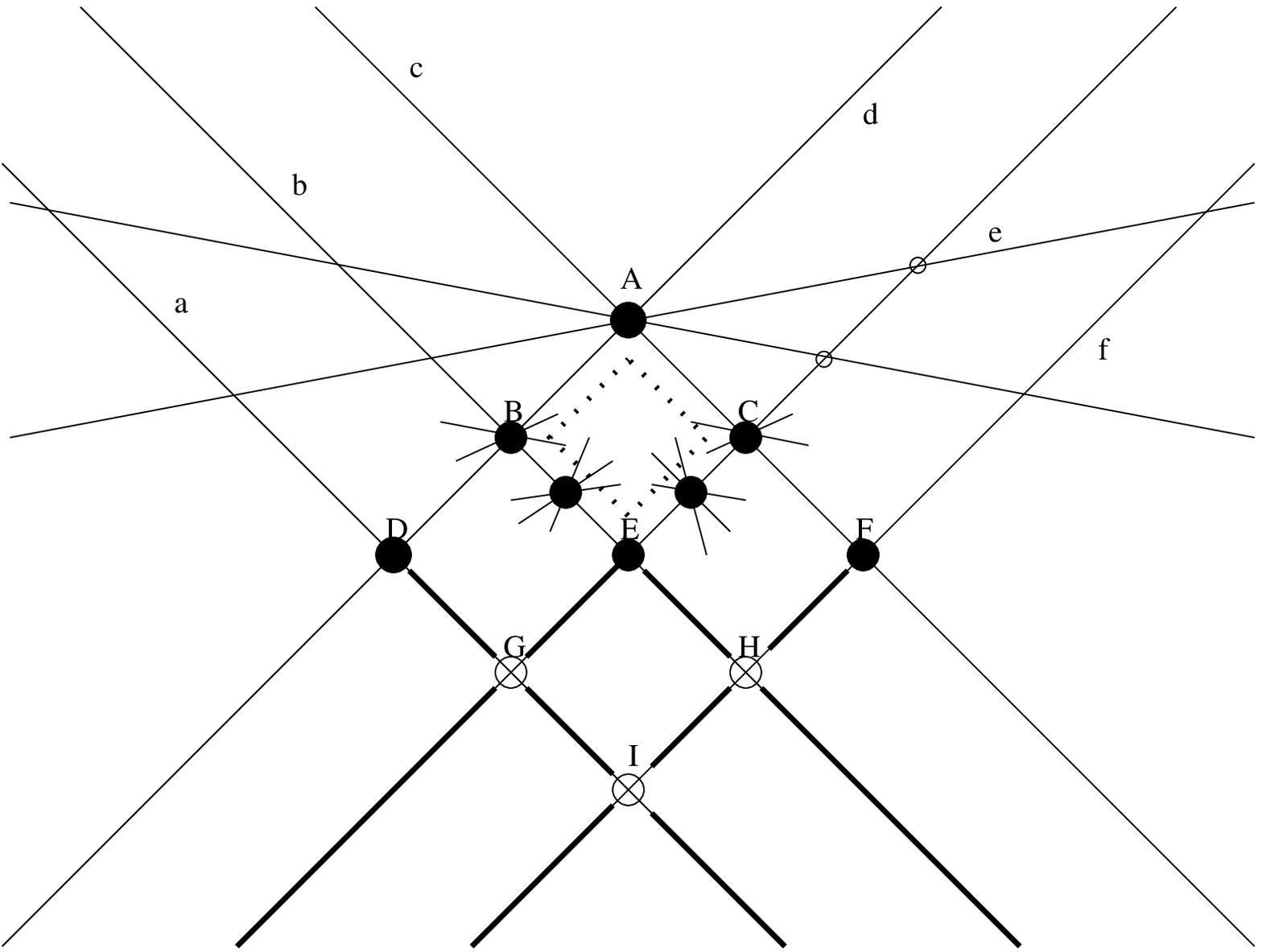}\label{fig:fall7b}}
  \caption{}
\end{figure}

\begin{lem}
  The situation, that the grid points $G$, $H$, and $I$ are
  $2$-crossings cannot occur. 
\end{lem}
\begin{proof}
  This is the situation of Figure~\ref{fig:fall7a}. We see that no
  $4$-crossing can lie on $a$ below $D$. No $4$-crossing can lie on
  $b$ below $E$. By symmetry no $4$-crossing can lie on $e$ below $E$,
  and no $4$-crossing can lie on $f$ below $F$. The drawing in
  Figure~\ref{fig:fall7b} shows that the $4$-crossings that are
  missing on $b$ and $e$ cannot both lie directly above $E$. We may
  assume that the missing $4$-crossing on $e$ lies above $C$; we call
  it $J$. The other line coming from $A$ crosses $e$ in a
  $2$-crossing, which means that only one further crossing lies on
  $e$. 
  
  Now we see that the five lines that come from $C$ and $J$ have to
  cross $f$ above $F$. They have the possibility to cross $f$ in the
  two missing $4$-crossings on $f$, which have to lie above $F$; we
  call them $K$ and $L$. This takes care of four of the five lines.
  However, one of the lines cannot go through $K$ or $L$, it yields a
  new $2$-crossing on $f$. So, all the crossings on $f$ are now
  determined.

  Now we take a closer look at $K$ and $L$. Of the two lines that come
  from $A$ at least one has to cross $f$ in $K$ or $L$ -- the other
  one could go through $J$ and leads to the above mentioned
  $2$-crossing. Together with the four lines that come from $C$ and
  $J$ five of the six lines that go through $K$ and $L$ are
  determined. So exactly one of the lines that go through $K$ and $L$
  does not cross $e$ in $C$ or $J$; we call this line $g$.
  
  To determine the place where $g$ crosses $e$, we first look at the
  lines that come from $F$. The two additional lines coming from $F$
  have to cross $e$ below $C$. However, we have already determined two
  $2$-crossings on $e$, and only one $4$-crossing lies below $C$. One
  of the lines leaving $F$ has to cross $e$ in $E$, and the other has
  to cross it in a $2$-crossing which either lies above $E$ or below
  $E$; we call these lines $h$ and $i$. If $i$ crosses $e$ above $E$,
  we call the resulting crossing $X$, and if $i$ crosses $e$ below
  $E$, we call the resulting crossing $Y$. All $2$-crossings are
  now determined, so $g$ has to cross $e$ in $E$.
  
  Now we have two cases. If $i$ crosses $e$ in $X$, $i$ has first to
  cross $g$. We call this crossing $M$. This crossing $M$, however,
  has to be a $2$-crossing. We have already determined all lines that
  enter the $2$-cell bordered by $C$, $E$, $F$, and $H$, and there are
  simply not enough of them to make $M$ a $4$-crossing. The
  $4$-crossings $E$ and $F$ are adjacent on $h$, no lines cross
  between them. And the segment between them borders two
  triangles which have a $2$-crossing as the remaining vertex. So, if
  we take $h$ to be the line at infinity, we are in the situation of
  Lemma~\ref{fig:fall2}, which settles this case.
  
  If $i$ crosses $e$ in $Y$, then $i$ has to cross $b$ in a
  $2$-crossing between $E$ and $H$, we call it $O$. Now we look at the
  points $M$, $F$, $O$, $E$, we are again in the situation of
  Lemma~\ref{fig:fall2} using $h$ as line at infinity.  
\end{proof}

So, with the proof of this lemma, we have settled all cases in which
six grid points were $4$-crossings. We already know, however that no
more than six grid points can be $4$-crossings
(Lemma~\ref{lem:amsix}), and we also know that no less than six grid
points can be $4$-crossings. This concludes the proof of
Theorem~\ref{thm:164}. 

\begin{rem}
  We can replace \emph{realizable} with \emph{pl-realizable} in
  Theorem~\ref{thm:164}.
\end{rem}

\section{Further Remarks}
\label{sec:further-results}

All $n_4$-configurations up to $n=17$ have been classified by Betten
and Betten \cite{MR1719478}. They have shown that there exist only
$19$ different $16_4$-configurations. This result gives us all
possible matroids that can lead to (pl-)realizable
$16_4$-configurations. Such a configuration can only be
pl-realizable, when the matroid is orientable. Using software, written
by the second author, we can decide whether a matroid is orientable.
The program tries to find a base orientation for the given matroid
that satisfy the Grassmann-Plücker-Relations. In all $19$ cases found
by Betten and Betten the matroid was not orientable, thus giving
another proof of Theorem~\ref{thm:164}.

We are optimistic that further arguments in connection with computer
support might lead to results in other cases as well.

\begin{figure}[h]
  \subfigure[$26_4$]{\igconf{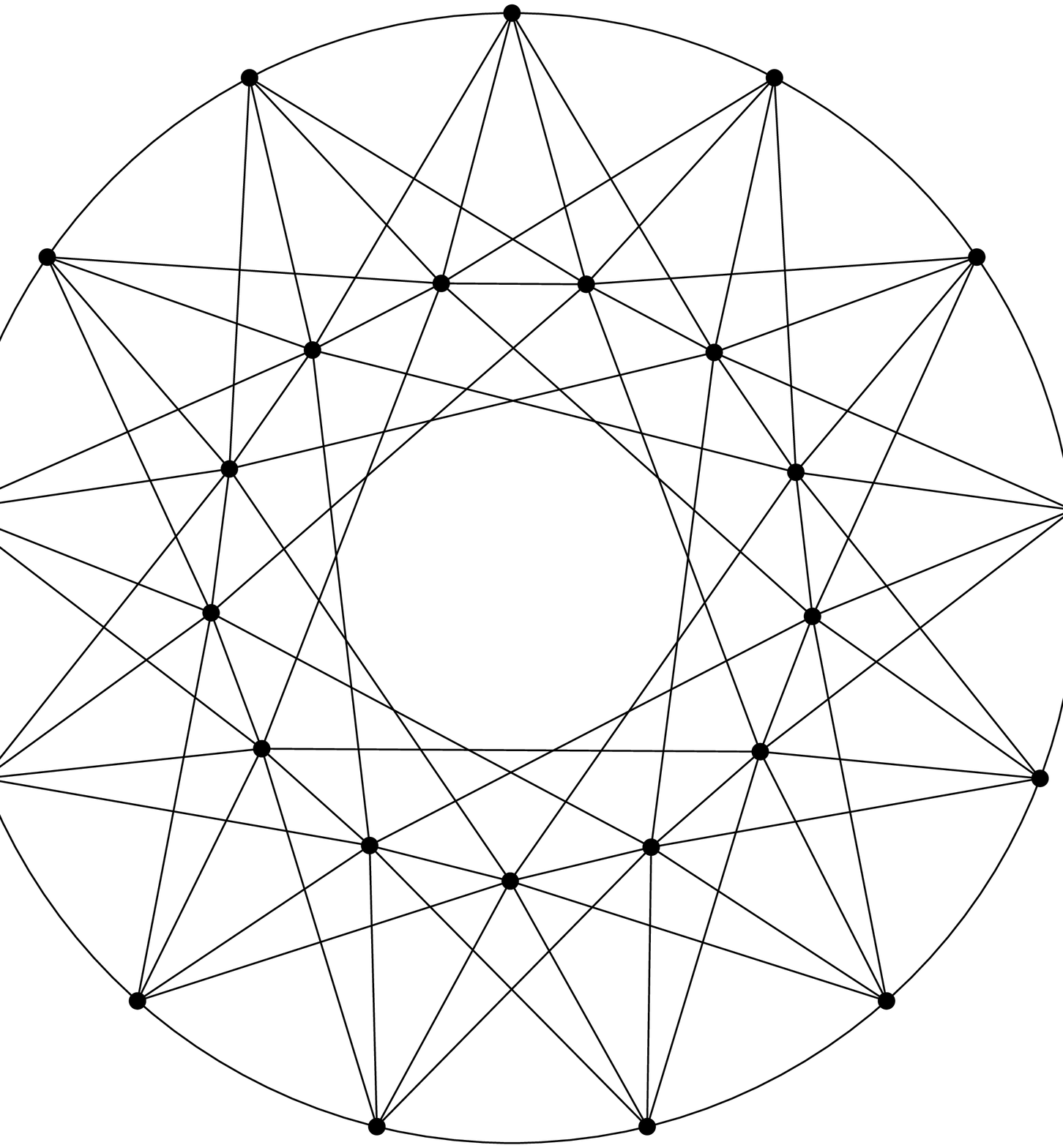}}\goodgap
  \subfigure[$32_4$]{\igconf{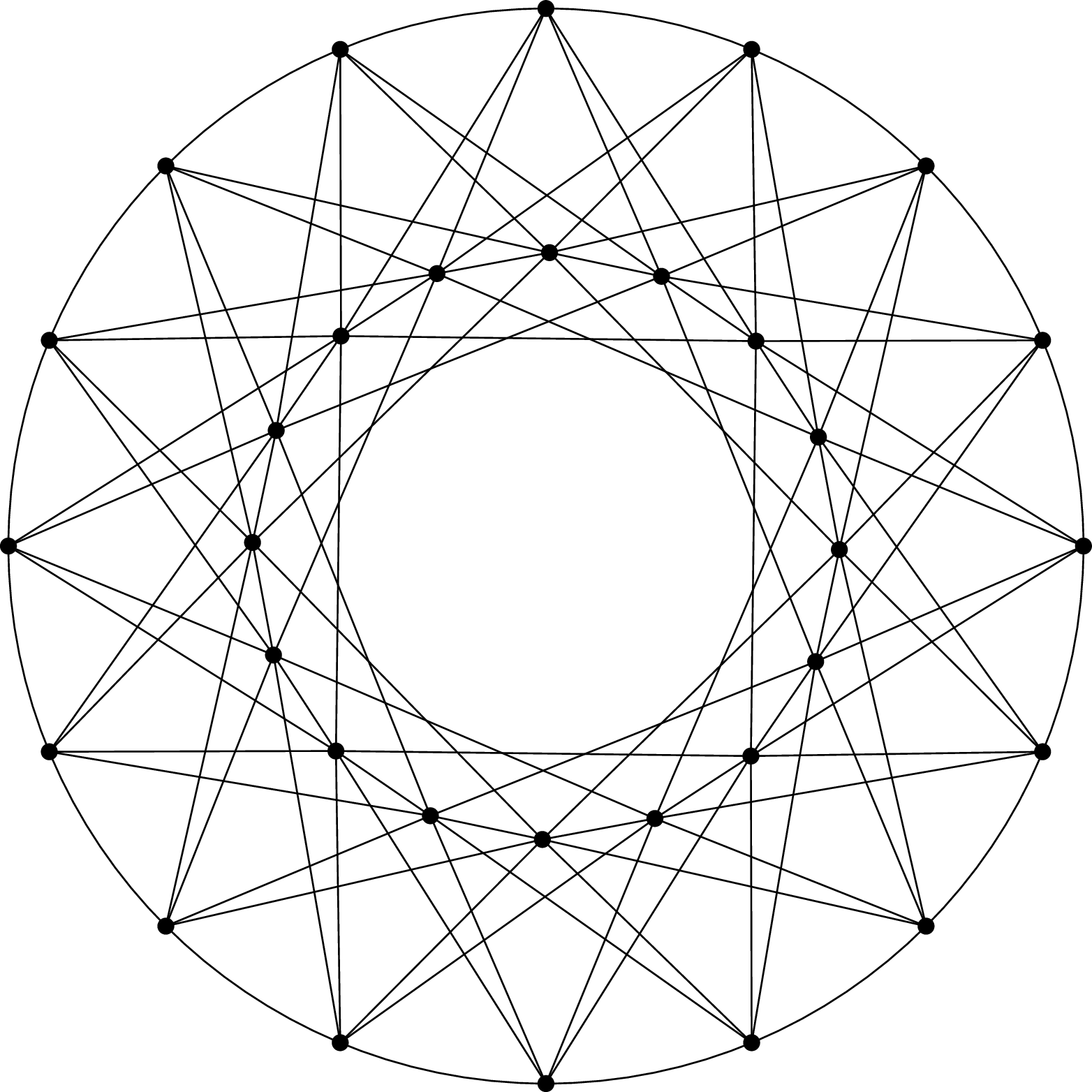}}\\
  \subfigure[$34_4$]{\igconf{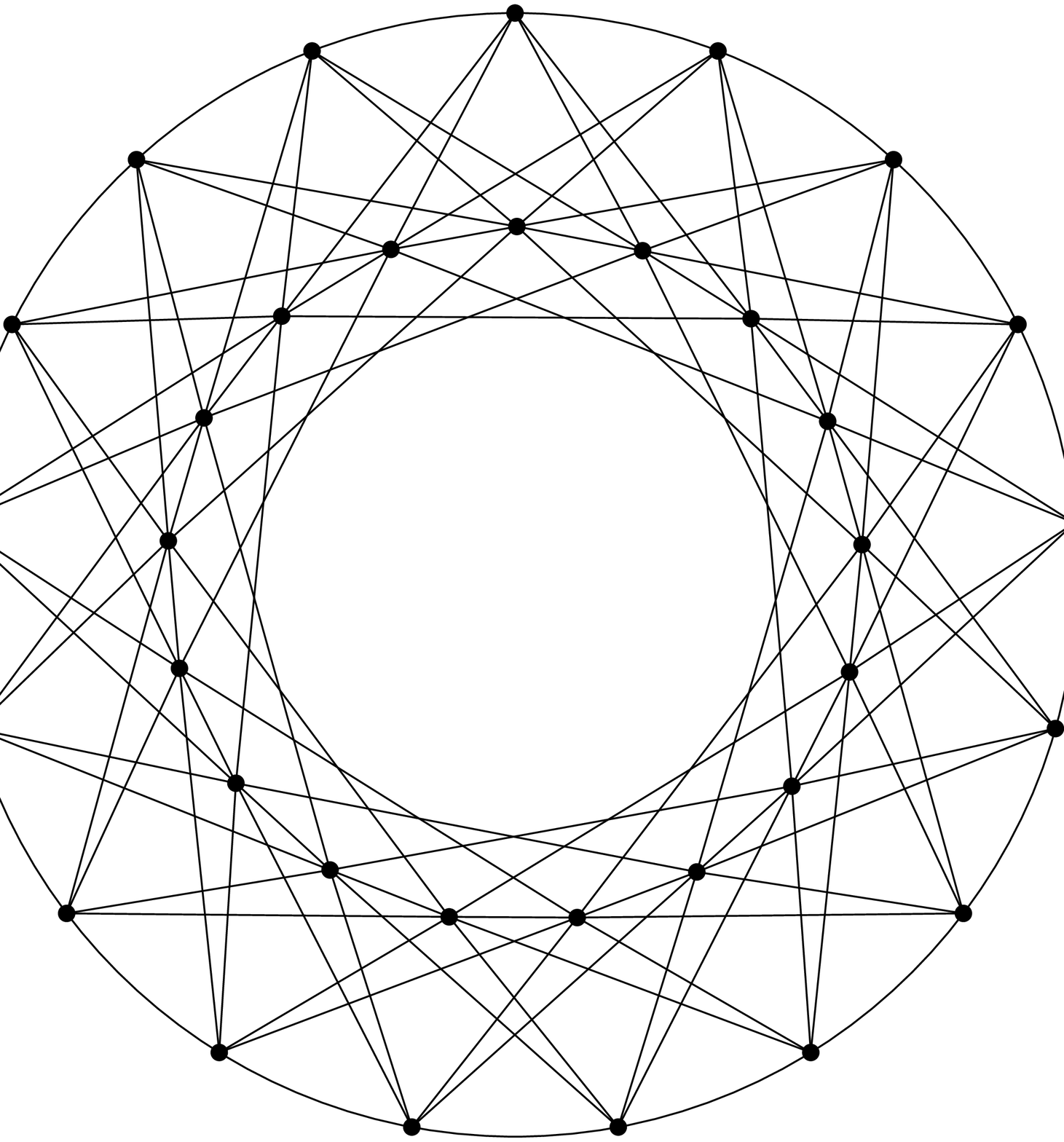}}\goodgap
  \subfigure[$38_4$]{\igconf{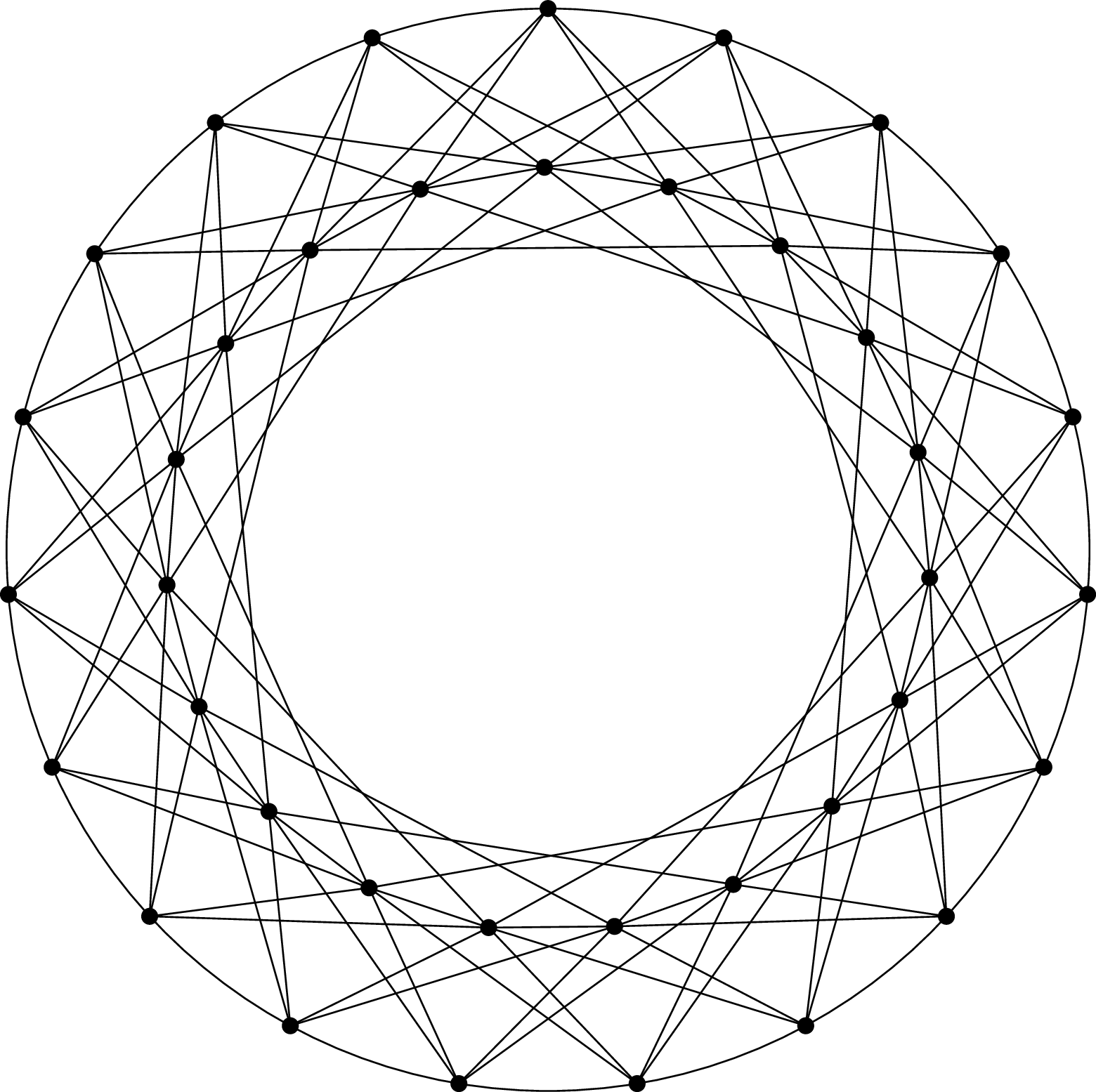}}
  \caption{Pseudoline Realizations by Branko Grünbaum}
  \label{fig:grun}
\end{figure}

After this article was written, Branko Gr{\"u}nbaum has sent
us his pseudoline arrangements of Figure~\ref{fig:grun}.

\newpage
\section*{Acknowledgements}
\label{sec:acknowledgements}

We thank Branko Grünbaum for his permission to publish his drawings of
pseudoline realizations for the cases $n=26,32,34,38$.

The second author was supported by a scholarship of the Deutsche
Telekom Foundation.


\begin{bibdiv}
\begin{biblist}
\bib{MR1719478}{article}{
    author={Betten, Anton},
    author={Betten, Dieter},
     title={Tactical decompositions and some configurations $v\sb 4$},
   journal={J. Geom.},
    volume={66},
      date={1999},
    number={1-2},
     pages={27\ndash 41},
      issn={0047-2468},
}
\bib{OM}{book}{
    author={Bj{\"o}rner, Anders},
    author={Las Vergnas, Michel},
    author={Sturmfels, Bernd},
    author={White, Neil},
    author={Ziegler, G{\"u}nter M.},
     title={Oriented Matroids},
    series={Encyclopedia of Mathematics and its Applications},
    volume={46},
   edition={2},
 publisher={Cambridge University Press},
     place={Cambridge},
      date={1999},
     pages={xii+548},
      isbn={0-521-77750-X},
}
\bib{CUP}{book}{
  author={Bokowski, J{\"u}rgen},
  title={Computational Oriented Matroids},
  publisher={Cambridge University Press},
  year={2005},
  note={to appear}
}
\bib{MR1845486}{article}{
    author={Bokowski, J{\"u}rgen},
    author={Mock, Susanne},
    author={Streinu, Ileana},
     title={On the Folkman-Lawrence topological representation theorem for
            oriented matroids of rank 3},
      note={Combinatorial geometries (Luminy, 1999)},
   journal={European J. Combin.},
    volume={22},
      date={2001},
    number={5},
     pages={601\ndash 615},
      issn={0195-6698},
}
\bib{MR1009366}{book}{
    author={Bokowski, J{\"u}rgen},
    author={Sturmfels, Bernd},
     title={Computational Synthetic Geometry},
    series={Lecture Notes in Mathematics},
    volume={1355},
 publisher={Springer-Verlag},
     place={Berlin},
      date={1989},
     pages={vi+168},
      isbn={3-540-50478-8},
}
\bib{MR0307027}{book}{
    author={Gr{\"u}nbaum, Branko},
     title={Arrangements and Spreads},
      note={Conference Board of the Mathematical Sciences Regional
            Conference Series in Mathematics, No. 10},
 publisher={American Mathematical Society Providence, R.I.},
      date={1972},
     pages={iv+114},

}
\bib{MR1913098}{article}{
    author={Gr{\"u}nbaum, Branko},
     title={Connected $(n\sb 4)$ configurations exist for almost all
            $n$---an update},
   journal={Geombinatorics},
    volume={12},
      date={2002},
    number={1},
     pages={15\ndash 23},
      issn={1065-7371},
}
\bib{MR1067273}{article}{
    author={Gr{\"u}nbaum, Branko},
    author={Rigby, J. F.},
     title={The real configuration $(21\sb 4)$},
   journal={J. London Math. Soc. (2)},
    volume={41},
      date={1990},
    number={2},
     pages={336\ndash 346},
      issn={0024-6107},
 }
\bib{leviPG}{article}{
author={Levi, Friedrich},
title={{Die Teilung der projektiven Ebene durch Gerade oder Pseudogerade}},
language={German},
journal={Sitz.ber. Sachs. Akad. Wiss. Leipz. Math.-Nat.wiss. Kl.},
volume={78},
pages={256-267},
year={1926},
}
\bib{leviCon}{book}{
author={Levi, Friedrich},
title={{Geometrische Konfigurationen. Mit einer Einführung in die Kombinatorische Flächentopologie}},
language={German},
publisher={{VIII${}+{}$310 S. 58 Abb. Leipzig, S. Hirzel. }},
year={1929},
}

 \bib{MR1217488}{book}{
    author={Orlik, Peter},
    author={Terao, Hiroaki},
     title={Arrangements of Hyperplanes},
    series={Grundlehren der Mathematischen Wissenschaften},
    volume={300},
 publisher={Springer-Verlag},
     place={Berlin},
      date={1992},
     pages={xviii+325},
     isbn={3-540-55259-6},
}
\end{biblist}
\end{bibdiv}
\end{document}